\documentclass[english,12pt]{article}
\usepackage{amsfonts}
\usepackage{amssymb}
\usepackage{amssymb, amsthm, amsmath}
\usepackage[latin1]{inputenc}
\usepackage{babel}
\usepackage{graphicx}
\usepackage{psfrag}
\usepackage{amstext}
\newtheorem{teorema}{\bf Theorem}

\newtheorem{observacao}{ \bf Remark}

\newtheorem{definicao}{ \bf Definition}

\newcommand{\abs}[1]{\lvert #1 \rvert}

\begin{document}
\title{ Minimal Surfaces in $S^3$ with Constant Contact Angle}
\author{Rodrigo Ristow Montes and Jose A. Verderesi
\thanks{ristow@mat.ufpb.br and javerd@ime.usp.br}}
\date{Departamento de Matem\'atica , \\
      Universidade Federal da Para\'iba, \\[2mm]
      BR-- 58.051-900 ~ Jo\~ao Pessoa, P.B., Brazil \\
and \\
Departamento de Matem\'atica Pura, \\
Instituto de Matem\'atica e Estat\'{\i}stica, \\
Universidade de S\~ao Paulo, \\
Caixa Postal 66281, \\[2mm]
BR--05315-970~ S\~ao Paulo, S.P., Brazil}

\maketitle
\renewcommand{\thefootnote}{\fnsymbol{footnote}}
\renewcommand{\thefootnote}{\arabic{footnote}}
\setcounter{footnote}{0}
\thispagestyle{empty}
\begin{abstract}
\noindent
  We provide a characterization of the Clifford Torus in $S^3$ via moving frames and contact structure equations. More precisely, we prove that minimal surfaces in $S^3$ with constant contact angle must be the Clifford Torus. Some applications of this result are then given, and some examples are discussed.
\end{abstract}
\smallskip

\noindent{\bf 2000 Math Subject Classification:} 53D10 - 53D35 - 53C40.
\section{Introduction}\mbox{}
The study of minimal surfaces played a formative role in the development of mathematics over the last two centuries. Today, minimal surfaces appear in various guises in diverse areas of mathematics, physics, chemistry and computer graphics, but have also been used in differential geometry to study  basic properties of immersed surfaces in contact manifolds. We mention for example the two papers, \cite{Blair} and \cite{MV}, where the authors classify Legendrian minimal surfaces in $S^5$ with constant Gaussian curvature. \\
Besides, interesting characterizations of the Clifford torus in spheres are given  in \cite{CI}, \cite{OP}, and \cite{TV}.\\
The scope of this note is to use a geometric invariant in order to study immersed surfaces in three dimensional sphere. This invariant (the contact angle $(\beta)$) is the complementary angle between the contact
distribution and the tangent space of the surface.\\
We show that the Gaussian curvature K of a minimal surface in $S^3$ with contact angle $\beta$ is given by:
\begin{eqnarray*}
K  &  = &  1 - \abs{\nabla\beta+e_1}^2 \nonumber
\end{eqnarray*}
Moreover, the contact angle satisfies the following Laplacian equation
\begin{eqnarray}
\Delta(\beta) &  =   & -\tan(\beta)\abs{\nabla\beta+2e_1}^2\label{eq:laplaciano1}\nonumber
\end{eqnarray}
where $e_1$ is the characteristic field defined in section 2 and  introduced by Bennequin, in \cite{A} pages 190 - 206.\\
Using the equations of Gauss and  Codazzi, we have proved the following two theorems:
\begin{teorema}\label{unicidade}
The Clifford Torus is the only minimal surface in $S^{3}$ with  constant contact angle.
\end{teorema}
\begin{teorema}\label{construcao}
The Clifford Torus is the only compact minimal  surface in $S^3$
with contact angle $0 \leq \beta < \frac{\pi}{2}$ ( or
$-\frac{\pi}{2} < \beta \leq 0$)
\end{teorema}
More in general, we have the following congruence result:
\begin{teorema}
Consider $S$ a Riemannian surface, \textbf{e} a vector field on $S$, and $\beta: S \rightarrow ]0,\frac{\pi}{2}[$ a function over $S$ that verifies the following equation:
\begin{eqnarray}
\Delta\beta = -\tan\beta (\abs{\nabla \beta}^2 + 4 (e(\beta)+1))
\end{eqnarray}
then there exist one and only one minimal immersion of $S$ into $S^3$ such that \textbf{e} is the characteristic vector fied, and $\beta$ is the contact angle of this immersion.
\end{teorema}
Finally, in section \ref{exemplos}, we give two examples of minimal surfaces in
$S^3$. The first one, we determine that the contact angle
($\beta$) of the Clifford Torus is ($\beta=0$) and the second one
we determine that the contact angle of the totally geodesic sphere
is ( $\beta=arc cos(x_2)$), and therefore, non constant.

\section{Contact Angle for Immersed Surfaces in  $S^{3}$}\label{sec:section2}
Consider in $\mathbb{C}^{2}$ the following objects:
\begin{itemize}
\item the Hermitian product: $(z,w)= z^1\bar{w}^1+z^2\bar{w}^2$;
\item the inner product: $\langle z,w \rangle = Re (z,w)$;
\item the unit sphere: $S^{3}=\big\{z\in\mathbb{C}^{2} | (z,z)=1\big\}$;
\item the \emph{Reeb} vector field in $S^{3}$, given by: $\xi(z)=iz$;
\item the contact distribution in $S^{3}$, which is orthogonal to $\xi$:
\[\delta_z=\big\{v\in T_zS^{3} | \langle \xi , v \rangle = 0\big\}.\]
\end{itemize}
Note that $\delta$ is invariant by the complex structure of $\mathbb{C}^{2}$.

Let now $S$ be an immersed orientable surface in $S^{3}$.
\begin{definicao}
The \emph{contact angle} $\beta$ is the complementary angle between the
contact distribution $\delta$ and the tangent space $TS$ of the
surface.
\end{definicao}
Let $(e_1,e_2)$ be  a local frame of $TS$, where $e_1\in
TS\cap\delta$. Then $\cos \beta = \langle \xi , e_2 \rangle
$.\\
Let  $(f_1=z^\bot$, $f_2=iz^\bot$ and $f_3=iz)$  be an orthonormal frame of $S^3$, where $z^\bot = (-\bar{z}_2, \bar{z}_1).$\\
The covariant derivative  is given by:
\begin{equation}
\begin{array}{ccc}\label{eq: deriv}
Df_1  &  =  &   \;w_1^2\,f_{2}+w^{2}\,f_{3}\\
Df_2  &  =  &   \;w_2^1\,f_{1}-w^{1}\,f_{3}\\
Df_3  &  =  & -w^2\,f_{1} +w^{1}\,f_{2}
\end{array}
\end{equation}
where $(w^{1},w^{2},w^{3})$ is the coframe associated to $(f_1,f_2,f_3)$.\\
Let $e_1$ be an  unitary vector  field in $TS\;\cap\;\delta$, where $\delta$ is the contact distribution.\\
Thus follows that:
\begin{equation}\label{eq: campos}
\begin{array}{ccl}
e_1  & = & f_{1} \\
e_2  & = & \sin(\beta)\,f_{2}+\cos(\beta)\,f_{3}\\
e_3  & = & -\cos(\beta)\,f_{2}+\sin(\beta)\,f_{3}
\end{array}
\end{equation}
where $\beta$ is the angle between $f_3$ and $e_2$,\,$(e_1,e_2)$ are tangent to $S$
and $e_3$ is normal to $S$
\section{Equations for the Gaussian Curvature and for the Laplacian of a Minimal Surface in $S^3$} \mbox{}
In this section, we will give formulas for the Laplacian and for
the Gaussian curvature of a minimal surface immersed in $S^3$.\\
The reader can see \cite{O}, and \cite{HG} for further details. \\
Let $(\theta^{1},\theta^{2},\theta^{3})$ be the coframe associated to $(e_1,e_2,e_3)$.\\
Thus, from equations (\ref{eq: campos}), it follows that:
\begin{equation}
\begin{array}{ccl}\label{eq: formas}
\theta^{1}    &   =   &   w_{1}\\
\theta^{2}    &   =   &   \sin(\beta)\,w_{2}+\cos(\beta)\,w_{3}\\
\theta^{3}    &   =   &   -\cos(\beta)\,w_{2}+\sin(\beta)\,w_{3}
\end{array}
\end{equation}
We know that $\theta^{3}=0$ on $S$, then we obtain the following equation:
\begin{equation}
\begin{array}{ccc}
sin(\beta)\,w^{3} & = & cos(\beta)w^{2}
\end{array}
\end{equation}
we have also
\begin{eqnarray}
w^2 & = & \sin \beta \theta^2 \nonumber\\
w^3 & = & \cos \beta \theta^2 \nonumber
\end{eqnarray}
It follows from (\ref{eq: formas}) that:
\begin{eqnarray}
d\theta^{1}     &   +     & \sin(\beta)(w_2^1-\cos(\beta)\theta^{2})\wedge\theta^{2}  =0\nonumber\\
d\theta^{2}     &   +     & \sin(\beta)(w_1^2+\cos(\beta)\theta^{2})\wedge\theta^{1} = 0\nonumber\\
d\theta^{3}     &   =     &   d\beta\wedge\theta^{2}-\cos(\beta)w_2^1\wedge\
w^{1}+(1+\sin^{2}(\beta))\theta^{1}\wedge\theta^{2} \nonumber
\end{eqnarray}
Therefore the connection form of $S$ is given by
\begin{eqnarray}
\theta_2^1     &   =     &   \sin(\beta)(w_2^1-\cos(\beta)\theta^{2})
\label{eq: formul}
\end{eqnarray}
Differentiating $e_3$ at the basis $(e_1,e_2)$, we have fundamental second forms coeficients
\begin{eqnarray}
De_3  &  =  &  \theta_3^1e_1+\theta_3^2e_2\nonumber
\end{eqnarray}
where
\begin{eqnarray}
\theta_3^1   &  =  & -\cos(\beta)w_2^1 - \sin^{2}(\beta)\theta^2 \nonumber\\
\theta_3^2   &  =  &  d\beta+\theta^1\nonumber
\end{eqnarray}
It follows from $d\theta^3=0$, that
\begin{eqnarray}
w_2^1(e_2)  &  =  &  -\frac {\beta_1}{\cos\beta}-\frac{(1+\sin^{2}\beta)}{\cos\beta}  \label{eq: formula1}
\end{eqnarray}
onde $d\beta(e_1)=\beta_1$.\\
The condition of minimality is equivalent to the following equation
\begin{eqnarray}
\theta_1^3\wedge\theta^2 - \theta_2^3\wedge\theta^1=0\nonumber
\end{eqnarray}
we have
\begin{eqnarray}
 w_2^1(e_1) & = \frac{\beta_2}{\cos(\beta)} \label{eq: formula2}
\end{eqnarray}
where $d\beta(e_2)=\beta_2$.\\
It follows from  (\ref{eq: formul}), (\ref{eq: formula1}) and (\ref{eq: formula2}),
\begin{eqnarray}
\theta_2^1  &  =  &  \tan(\beta)(\beta_2\theta^{1}-(\beta_{1}+2)\theta^{2})\nonumber \\
\theta_3^1  &  =  &  -\beta_2\theta^1+(\beta_1+1)\theta^2\nonumber\\
\theta_3^2  &  =  &   (\beta_1+1)\theta^1+\beta_2\theta^2\nonumber
\end{eqnarray}
If $J$ is the complex structure of $S$ we have $Je_1=e_2$ e $Je_2=-e_1$.\\
Using $J$, the forms above reduce to:
\begin{eqnarray}
\theta_2^1  &  =  &  \tan\beta(d\beta\circ J-2\theta^2)\nonumber \\
\theta_3^1  &  =  &  -d\beta\circ J+\theta^2\\\nonumber
\theta_3^2  &  =  &   d\beta+\theta^1\nonumber
\end{eqnarray}
Gauss equation is
\begin{eqnarray}
d\theta_1^2=\theta^2\wedge\theta^1+\theta_1^3\wedge\theta_2^3 \nonumber
\end{eqnarray}
which implies
\begin{equation}
\begin{array}{lcl}
d\theta_2^1   &   =  &  (|\nabla\beta|^2+2\beta_1)\;(\theta^{2}\wedge\theta^{1}) \label{eq: curva2}
\end{array}
\end{equation}
and therefore
\begin{eqnarray}
K  &  = &  1 - |\nabla\beta+e_1|^2 \nonumber
\end{eqnarray}
Differentiating $\theta_2^1$, we have
\begin{equation}
\begin{array}{lcl}
d\theta_2^1  &  =  & \quad \sec^{2}(\beta)(|\nabla\beta|^2+2\beta_1)(\theta^{2}\wedge\theta^{1})\\
             &     & + (\tan(\beta)\Delta(\beta)+2\tan^{2}(\beta)(\beta_1+2))(\theta^{2}\wedge\theta^{1}) \label{eq: curva1}
\end{array}
\end{equation}
Using  (\ref{eq: curva2}) and (\ref{eq: curva1}), we obtain the following formula for the Laplacian of $S$
\begin{equation}
\begin{array}{lcl}
\Delta(\beta) &  =   & -\tan(\beta)((\beta_1+2)^2+\beta^2_2)\label{eq:laplaciano1}
\end{array}
\end{equation}
Or
\begin{eqnarray}
\Delta(\beta) &  =   & -\tan(\beta)|\nabla\beta+2e_1|^2 \nonumber
\end{eqnarray}
Codazzi equations are
\begin{eqnarray}\label{eq: codazzi 1}
d\theta_1^3+\theta_2^3{\wedge}\theta_1^2=0\nonumber\\
\label{eq: codazzi 2}
d\theta_2^3+\theta_1^3{\wedge}\theta_2^1=0\nonumber
\end{eqnarray}
A straightforward computation in the first equation gives
(\ref{eq:laplaciano1}) and the second equation is always verified.
\section{Main Results}
\subsection{Proof of the Theorem 1} \mbox{}
Suppose that  $\beta$ is constant, it follows from  (\ref{eq: curva2}) that
$d\theta_1^2=0$ and, therefore, $K=0$,ie., Gaussian curvature is identically
null, hence $S \subset S^3$ is the Clifford Torus, which prove the Theorem 1.
\subsection{Proof of the Theorem 2} \mbox{}
For  $0 \leq \beta < \frac{\pi}{2}$, we have $\tan \beta \geq  0$, hence
$\Delta(\beta) \leq 0$ and using that  $S$ is a compact surface, we conclude
by Hopf's Lemma that $\beta$ is constant, and therefore, $K=0$ and $S$ is the
Clifford Torus, which prove the Theorem 2.
\subsection{Proof of the Theorem 3}
Let $S$ be an orientable surface in $S^3$, and let \textbf{e} be an unit vector field on $S$. \\
We choose an orthonormal positive basis $(e_1,e_2)$ with $e_1=e$, and let $(\theta^1,\theta^2)$ be a coframe on $S$.\\
For each function $\beta : S \rightarrow ]0,\frac{\pi}{2}[$ that satisfies the following Laplacian equation:
\begin{eqnarray}
\Delta(\beta) &  =   & -\tan(\beta)\abs{\nabla\beta+2e_1}^2\nonumber
\end{eqnarray}
We define the following  fundamental  second form:
\begin{eqnarray}
\theta_1^3 & = & (d\beta + \theta^1)\circ J\nonumber \\
\theta_2^3 & = & \phantom{-} -(d\beta + \theta^1)\nonumber \\
\end{eqnarray}
Now, the proof follows from Gauss-Codazzi equations.
\section{ Examples}\label{exemplos}
Examples of minimal surfaces in $S^3$ was discovered by Lawson, in \cite{BL}.\\
Here we will use the notion of the contact angle to give a characterization of these known examples.
\subsection{Contact Angle of Clifford Torus in $S^3$}\mbox{}\\
Let us consider the following torus in $S^3$:
\begin{eqnarray*}
T^2=\lbrace(z_1,z_2)\in{C^2}/z_1\bar{z_1}=\frac{1}{2},z_2\bar{z_2}=\frac{1}{2}\rbrace
\end{eqnarray*}
Let $f$ be the following immersion:
\begin{eqnarray*}
f(u_1,u_2)=\frac{\sqrt 2}{2}(e^{iu_1},e^{iu_2})
\end{eqnarray*}
Tangent space $T(T^2)$ is  given by $\frac{\partial}{\partial u_1}$ and
$\frac{\partial}{\partial u_2}$, thus we have:
\begin{eqnarray*}
a\frac{\partial}{\partial u_1}+b\frac{\partial}{\partial u_2}={\lambda}z^{\bot}
\end{eqnarray*}
using the condition above and the fact that $| \lambda |=1$, we obtain:
\begin{eqnarray*}
\lambda   &   =    &   ie^{i(u_1+u_2)}
\end{eqnarray*}
Unit vector fields are:
$$\left\{
\begin{array}{lll}
e_1=ie^{i(u_1+u_2)}z^\bot\\
e_2=iz\\
e_3=iz^{\bot}
\end{array}\right.$$
The contact angle is the angle between $e_2$ and $f_3$  ,
$$
\begin{array}{ccl}
cos(\beta) & = & \langle e_2,f_3 \rangle\\
           & = &  1
\end{array}
$$
Therefore, the contact angle is:
$$
\beta=0
$$
Fundamental second form is given by the following:\\
\begin{center}
$A=\left[\begin{array}{cc}
0 & -1\\ -1  & 0\\
\end{array}\right]$
\end{center}
\subsection{Minimal surface in $S^3$ with non constant contact angle}\mbox{}
Let us consider the following surface in $S^3$:\\
$$\left\{
\begin{array}{lcc}
z_2-\overline{z}_2              & =  &   0\\
(x_1)^2+(y_1)^2+(x_2)^2+(y_2)^2 & =  &   1
\end{array}\right.$$\\
We see that the vector fields are:
$$\left\{
\begin{array}{lll}
e_1=\frac{1}{\sqrt{1-x_2^{2}}}(-x_1x_2,-y_1x_2,1-x_2^{2},0)\\
e_2=\frac{1}{\sqrt{x_1^{2}+y_1^{2}}}(y_1,-x_1,0,0)\\
e_3=(0,0,0,1)
\end{array}\right.$$
The contact angle is the angle between $e_2$ and $f_3$, that is,
$$
\begin{array}{ccl}
cos(\beta) & = & \langle e_2,f_3 \rangle\\
           & = &  x_2
\end{array}
$$
Therefore, the contact angle is non constant:
$$
\beta= arc \, cos\, (x_2)
$$
\begin{observacao}
For higher dimensions, when we have a compact minimal surface immersed in $S^5$, we proved, in \cite{MV}, that the case $\beta = \frac{\pi}{2}$  gives an alternative proof of the classification of a Theorem from 
Blair in \cite{Blair}, for Legendrian minimal surfaces in $S^5$ with
constant Gaussian curvature
\end{observacao}


\end{document}